\documentclass[12pt,reqno]{amsart}
\usepackage{amssymb,amsmath,amsthm,longtable,multirow,array}

\usepackage{color}
\usepackage{mathrsfs}
\usepackage{tikz}
\usepackage{subeqnarray}
\usepackage{cases}

\numberwithin{equation}{section}

\newtheorem{Theorem}{Theorem}[section]
\newtheorem{Proposition}[Theorem]{Proposition}

\newtheorem{Lemma}[Theorem]{Lemma}
\newtheorem{Corollary}[Theorem]{Corollary}
\theoremstyle{definition}
\newtheorem{Definition}[Theorem]{Definition}
\theoremstyle{remark}
\newtheorem{Example}[Theorem]{Example}

\newtheorem*{remark*}{Remark}

 \providecommand\CC{{\mathbb C}}

\providecommand\brs{\begin{remark*}}
\providecommand\ers{\end{remark*}}

\providecommand\be{\begin{enumerate}}
\providecommand\ee{\end{enumerate}}
\providecommand\bT{\begin{Theorem}}
\providecommand\eT{\end{Theorem}}
\providecommand\bP{\begin{Proposition}}
\providecommand\eP{\end{Proposition}}
\providecommand\bD{\begin{Definition}}
\providecommand\eD{\end{Definition}}
\providecommand\bE{\begin{Example}}
\providecommand\eE{\end{Example}}
 \providecommand\bL{\begin{Lemma}}
\providecommand\eL{\end{Lemma}}
\providecommand\bC{\begin{Corollary}}
\providecommand\eC{\end{Corollary}}

\providecommand\bpp{\begin{proof}} \providecommand\epp{\end{proof}}
\providecommand\bee{\begin{equation}}
\providecommand\eee{\end{equation}}

  \providecommand\beqq{\begin{eqnarray*}}
\providecommand\eeqq{\end{eqnarray*}}
\providecommand\a{\alpha}

 \providecommand\rt{\rightarrow}
\providecommand\?{\infty}

\providecommand\besp{\begin{split}}
\providecommand\eesp{\end{split}}
\providecommand\bay{\begin{array}}
\providecommand\eay{\end{array}}

\providecommand\bes{\begin{equation}\begin{split}}
\def\a{\alpha}
\def\b{\beta}

\def\f#1#2{\frac{#1}{#2}}

\def\t{\tau}

\def\<{\leq}
\begin{document}

\title[Orthogonal Trigonometric Polynomial]{The Strong Asymptotic Analysis of the first kind Orthogonal Trigonometric Polynomial}
\author{Han Huili}
\author{Hua Liu}
\author{Wang Yufeng}

\begin{abstract} In this paper we study the asymptotic analysis of the orthogonal trigonometric polynomials by the  Riemann-Hilbert problem for the periodic analytic functions.

\end{abstract}

\address{Han Huili, Ningxia University}
\address{Liu Huadaliuhua@163.com, Department of Mathematics, Tianjin
University of Technology and Education, Tianjin 300222, China}
\address{Wang Yufeng, Wuhan University}

\maketitle

\bigskip

\section{orthogonal trigonometric polynomial}
\vskip5mm
If a nonnegative locally-integrable function $w(x)$ defined on the real axis $\mathbb{R}$ satisfying
\begin{equation}
w(x+2\pi)=w(x)\,\,\,\mathrm{for}\,\,x\in \mathbb{R}\,\,\,\mathrm{and}\,\,\,\int_0^{2\pi}w(x)\mathrm{d}x>0,
\end{equation}
such a function $w$ is called a $2\pi-$periodic weight. And the real inner product is defined by
\begin{equation}
\langle f,g\rangle=\int_0^{2\pi}f(x)g(x)w(x)\mathrm{d}x,
\end{equation}
which induces the norm
\begin{equation}
\|f\|_2=\int_0^{2\pi}|f(x)|^2w(x)\mathrm{d}x.
\end{equation}

By use of the inner product (1.2), the Gram-Schmidt orthogonalization of the following ordered trigonometric monomials
\begin{equation}
  1,\cos t,\sin t,\cdots,\cos nt,\sin nt,\cdots
\end{equation}
leads to the system of orthonormal trigonometric polynomials $\{\omega_n,n=0,1,2,\cdots\}$, where
\begin{equation}
\omega_{2n}(x)=\a_n\cos nx+\displaystyle\sum^{n-1}_{k=1}(a_k \cos kx+b_k\sin kx)+a_0\,\,\,\mathrm{with}\,\,\,\a_n>0
\end{equation}
and
\begin{equation}
\omega_{2n+1}(x)=\b_n\sin nx+a_n\cos nx+\displaystyle\sum^{n-1}_{k=1}(a_k \cos kx+b_k\sin kx)+a_0\,\,\,\mathrm{with}\,\,\,\b_n>0.
\end{equation}
The trigonometric polynomial $\omega_{2n}$ defined by (1.5) is usually called the first kind orthogonal trigonometric polynomial (OTP) and $\a_n$ is said to be its leading coefficient.
And $\omega_{2n+1}$ defined by (1.6) is similarly called the second kind OTP and $\b_n$ is the corresponding leading coefficient.

For convenience, we define two real vector spaces according to the number of the base
\begin{equation}
T_{2n}(\mathbb{R})=\left\{a_n\cos nx+\sum^{n-1}_{k=0}(a_k \cos kx+b_k\sin kx)+a_0,\,a_j,b_j\in \mathbb{R}\right\},
\end{equation}
\begin{equation}
T_{2n+1}(\mathbb{R})=\left\{\sum^{n}_{k=0}(a_k \cos kx+b_k\sin kx)+a_0,\,a_j,b_j\in \mathbb{R}\right\}.
\end{equation}

Therefore, one has
 \begin{equation}
 \left\{
 \begin{array}{l}
\displaystyle\left\langle\omega_{2n},\cos kx\right\rangle=0,\,\,\, k=0,1,2\cdots,n-1,\\[3mm]
\displaystyle\left\langle\omega_{2n},\sin kx\right\rangle=0,\,\,\, k=1,2\cdots,n-1,\\[3mm]
\displaystyle\a^2_n=\left\langle \omega_{2n},\omega_{2n}\right\rangle=1,
  \end{array}\right.
\end{equation}
where $\omega_{2n}$ is the first kind OTP defined by (1.5). Conversely, if $\omega_{2n}$ satisfies (1.9), one easily knows that it is just the first kind OTP defined by (1.5). Further, the first two equalities in (1.9)
are equivalent to
\begin{equation}
\omega_{2n}\bot T_{2n-1}(\mathbb{R}),
\end{equation}
where $ T_{2n-1}(\mathbb{R})$ is defined by (1.8).

Further, we define the first kind  monic OTP
\begin{equation}
\varpi_{2n}(x)=\frac{\omega_{2n}(x)}{\a_n},
\end{equation}
where $\a_n$ is the corresponding leading coefficient.

A $2\pi-$periodic weight $w$ is called strictly-positive analytic periodic weight if $w\in A(\mathbb{R})$ and $w(x)>0,\,x\in \mathbb{R}$. In this section, we always assume that
$w$ is strictly-positive analytic periodic weight. Thus, there exists $\rho>0$ such that
\begin{equation}
 \left\{
 \begin{array}{l}
\displaystyle w\in A(\square_{-\rho,\rho}),\\[3mm]
\displaystyle w(z)\neq 0,\,\,\,z\in \square_{-\rho,\rho},
  \end{array}\right.
\end{equation}
where $A(\square_{-\rho,\rho})$ is the set of analytic functions on the rectangle domain $\square_{-\rho,\rho}$ defined by
\begin{equation}\square_{a,b}=\big\{z,\mathrm{Re}z\in (0.2\pi),\,\mathrm{Im}z\in(a,b)\big\}\,\,\,\mathrm{for}\,\,\,a<b.\end{equation}

\section{Characterization of the first kind OTP}

In this section, we will give the Fokas-Its-Kitaev characterization of the first kind OTP.
Now, we state the homogeneous periodic Riemann-Hilbert problem: Find a sectionally-analytic $2\pi-$periodic function ${\bf Y}(z)$ satisfying the following conditions
\begin{equation}
\left\{
\begin{array}{ll}
{\bf Y}^+(x)={\bf Y}^-(x)\left(\!\!
\begin{array}{cc}
1&e^{-inx}w(x)\\
0&1
\end{array}
\!\!\right),\,\,\,\,&x\in [0,2\pi],
\\[5mm]
{\bf Y}(z){\bf \Xi_1}(z)\rightarrow\mathbf{I},\,\,\,&z\rightarrow +\infty i,\\[3mm]
{\bf Y}(z){\bf \Xi_2}(z)\rightarrow\mathbf{I},\,\,\,&z\rightarrow -\infty i,
\end{array}\right.
\end{equation}
where $\mathbf{I}$ is the $2\times 2$ identity matrix and
\begin{equation}
\mathbf{\Xi}_1(z)=\left(
            \begin{array}{cc}
              \f1{\cos nz} & 0 \\[3mm]
              0 & e^{iz} \\
            \end{array}  \right),\,\,\,
  \mathbf{\Xi}_2(z)= \left(
            \begin{array}{cc}
              \f1{\cos nz} & 0 \\[3mm]
              0 & e^{i(2n-1)z} \\
            \end{array} \right).
\end{equation}

Let
\begin{equation}
{\bf Y}(z)=\left(
          \begin{array}{cc}
            Y_{1,1}(z) & Y_{1,2}(z) \\
            Y_{2,1}(z) & Y_{2,2}(z) \\
          \end{array}
        \right),\,\,\,\,z\in \mathbb{C}\setminus\mathbb{R},
\end{equation}
and the periodic matrix Riemann-Hilbert problem (2.1) is equivalent to the system of the following four scalar Riemann-Hilbert problems
\begin{equation}
\left\{
\begin{array}{ll}
 Y_{1,1}^+(x)= Y_{1,1}^-(x),\,\,\,\,&x\in [0,2\pi],
\\[3mm]
\displaystyle\frac{Y_{1,1}(z)}{\cos nz}\rightarrow 1,\,\,\,&z\rightarrow +\infty i,\\[3mm]
\displaystyle\frac{Y_{1,1}(z)}{\cos nz}\rightarrow 1,\,\,\,&z\rightarrow -\infty i,
\end{array}\right.
\end{equation}
\begin{equation}
\left\{
\begin{array}{ll}
 Y_{1,2}^+(x)= Y_{1,2}^-(x)+e^{-inx}w(x) Y_{1,1}^-(x),\,\,\,\,&x\in [0,2\pi],
\\[3mm]
\displaystyle e^{iz}Y_{1,2}(z)\rightarrow 0,\,\,\,&z\rightarrow +\infty i,\\[3mm]
\displaystyle e^{i(2n-1)z}Y_{1,2}(z)\rightarrow 0,\,\,\,&z\rightarrow -\infty i,
\end{array}\right.
\end{equation}
\begin{equation}
\left\{
\begin{array}{ll}
 Y_{2,1}^+(x)= Y_{2,1}^-(x),\,\,\,\,&x\in [0,2\pi],
\\[3mm]
\displaystyle\frac{Y_{2,1}(z)}{\cos nz}\rightarrow 0,\,\,\,&z\rightarrow +\infty i,\\[3mm]
\displaystyle\frac{Y_{2,1}(z)}{\cos nz}\rightarrow 0,\,\,\,&z\rightarrow -\infty i,
\end{array}\right.
\end{equation}
and
\begin{equation}
\left\{
\begin{array}{ll}
 Y_{2,2}^+(x)= Y_{2,2}^-(x)+e^{-inx}w(x) Y_{2,1}^-(x),\,\,\,\,&x\in [0,2\pi],
\\[3mm]
\displaystyle e^{iz}Y_{2,2}(z)\rightarrow 1,\,\,\,&z\rightarrow +\infty i,\\[3mm]
\displaystyle e^{i(2n-1)z}Y_{2,2}(z)\rightarrow 1,\,\,\,&z\rightarrow -\infty i.
\end{array}\right.
\end{equation}

In order to deal with these periodic Riemann-Hilbert problems, we need introduce the following periodic Cauchy-type integral operator\cite{W-W}
\begin{equation}
\mathcal{C}_w[f](z)=\f1{4\pi i}\int_0^{2\pi}f(t)\cot\f{t-z}2w(t)\mathrm{d}t,\,\,\,z\in \mathbb{C}\setminus\mathbb{R}
\end{equation}
with $f\in H([0,2\pi])$\cite{lu}. For convenience, we similarly define two complex vector spaces of trigonometric polynomials
\begin{equation}
T_{2n}(\mathbb{C})=\left\{a_n\cos nx+\sum^{n-1}_{k=0}(a_k \cos kx+b_k\sin kx)+a_0,\,a_j,b_j\in \mathbb{C}\right\},
\end{equation}
\begin{equation}
T_{2n+1}(\mathbb{C})=\left\{\sum^{n}_{k=0}(a_k \cos kx+b_k\sin kx)+a_0,\,a_j,b_j\in \mathbb{C}\right\}.
\end{equation}
Clearly, $T_{2n}(\mathbb{R})\subset T_{2n}(\mathbb{C})$ and $T_{2n+1}(\mathbb{R})\subset T_{2n+1}(\mathbb{C})$.
\vspace{3mm}

\noindent{\bf Theorem 2.1}\,\,\,
 {\it The matrix Riemann-Hilbert
problem $(2.1)$ has the unique solution expressed by
\begin{equation}\label{3.18}
\mathbf{Y}(z)\!=\!\!\left(
\!\!\begin{array}{cc}
\varpi_{2n}(z)&\displaystyle e^{-inz}\mathcal{C}_w[\varpi_{2n}](z)\\[3mm]
a_n\varpi_{2n-1}(z)&\displaystyle
a_ne^{-inz}\mathcal{C}_w[\varpi_{2n-1}](z)
\end{array}\!\!\right)\!,\,\,z\!\in\! \mathbb{C}\setminus \mathbb{R}
\end{equation}
with
\begin{equation}
a_n=\frac{2\pi i}{\beta_{n-1}^2},
 \end{equation}
where $\varpi_{2n},\,\varpi_{2n-1}$ are the monic OLPs, and $\beta_{n-1}$ is the leading coefficients of OTPs defined by (1.6).
}\vspace{3mm}

{\bf Proof:}\,\,\,First, we will solve scalar Riemann-Hilbert problems (2.4) and (2.5). By Theorem 4.1 in \cite{W-W}, the solution of Riemann-Hilbert problems (2.4)
can be expressed by
\begin{equation}
  Y_{1,1}(z)=\sum^{n}_{k=0}(a_k\cos kz+b_k\sin kz)
\end{equation}
with $a_j,\,b_j\in \mathbb{C}$ for $j=0,1,\cdots,n$. This leads to
\begin{equation}\lim_{z\rt\pm\?i}\f{Y_{1,1}(z)}{\cos nz}=a_n\pm b_ni,\end{equation}
by the growth conditions in (2.4). Therefore, one has
\begin{equation}
  Y_{1,1}(z)=t_{2n}(z)\in T_{2n}(\mathbb{C}),
\end{equation}
and the leading coefficient of $t_{2n}(z)$ is $1$.

{Let $\widetilde{Y}_{1,2}(z)$ satisfy }
\begin{equation}
\left\{
\begin{array}{l}
 \widetilde{Y}_{1,2}^+(x)= \widetilde{Y}_{1,2}^-(x)+t_{2n}(x)w(x),\,\,\,\,x\in [0,2\pi],
\\[3mm]
\displaystyle \widetilde{Y}_{1,2}(\pm \infty i)=0.
\end{array}\right.
\end{equation}
By Theorem 4.1 in \cite{W-W},
\begin{equation}
\widetilde{Y}_{1,2}(z)=\mathcal{C}_w[t_{2n}](z)=\f1{4\pi i}\int_0^{2\pi}t_{2n}(\tau)\cot\f{\tau-z}2w(\tau)\mathrm{d}\tau,\,\,\,z\in \mathbb{C}\setminus\mathbb{R}
\end{equation}
is the unique solution of (2.16)
if it is solvable. Therefore, if (2.5) is solvable, one has
\begin{equation}
Y_{1,2}(z)=e^{-inz}\widetilde{Y}_{1,2}(z)=\f{e^{-inz}}{4\pi i}\int_0^{2\pi}t_n(\tau)\cot\f{\tau-z}2w(\tau)\mathrm{d}\tau,\,\,\,z\in \mathbb{C}\setminus\mathbb{R}.
\end{equation}

Observe
\begin{equation}
\cot\f{\tau-z}2=\left\{
\begin{array}{ll}
 \displaystyle i\left(1+2\sum^\?_{k=1}e^{ikz}e^{-i\tau k}\right),\,\,\,\,&\mathrm{Im} z>0
\\[3mm]
\displaystyle -i\left(1+2\sum^\?_{k=1}e^{-ikz}e^{i\tau k}\right),\,\,\,\,&\mathrm{Im} z<0
\end{array}\right.\,\,\,\mathrm{for}\,\,\,\tau\in [0,2\pi].
\end{equation}
Now, inserting (2.19) into (2.18), one has the the following Fourier
expansion
\begin{equation}
Y_{1,2}(z)=\left\{
\begin{array}{l}
 \displaystyle  \f1{4\pi}\!\!\int^{2\pi}_0t_{2n}(\tau)w(\tau)\mathrm{d}\tau e^{-inz}+\sum^\?_{k=1}\f1{2\pi}\!\!\int^{2\pi}_0t_{2n}(\tau)e^{-ik\tau}w(\tau)\mathrm{d}\tau e^{-i(n-k)z},\,\,\,\,\mathrm{Im} z>0,
\\[3mm]
\displaystyle -\f1{4\pi}\!\!\int^{2\pi}_0t_{2n}(\tau)w(\tau)\mathrm{d}\tau e^{-inz}-\sum^\?_{k=1}\f1{2\pi}\!\!\int^{2\pi}_0t_{2n}(\tau)e^{ik\tau}w(\tau)\mathrm{d}\tau e^{-i(n+k)z},\,\,\,\,\mathrm{Im} z<0.
\end{array}\right.
\end{equation}
Combining the growth conditions in (2.5) with (2.20), one easily knows, if and only if
\begin{equation}
\left\{
\begin{array}{l}
 \displaystyle \int^{2\pi}_0t_{2n}(\tau)e^{-ik\tau}w(\tau)\mathrm{d}\tau=0,
\\[3mm]
\displaystyle  \int^{2\pi}_0t_{2n}(\tau)e^{ik\tau}w(\tau)\mathrm{d}\tau=0,
\end{array}\right.\ k=0,1,\cdots,n-1,
\end{equation}
the unique solution of (2.5) can be expressed by (2.18) or (2.20). Obviously, (2.21) is equivalent to
\begin{equation}
 \left\{
 \begin{array}{l}
\displaystyle\left\langle t_{2n},\cos k\tau\right\rangle=0,\,\,\, \\[3mm]
\displaystyle\left\langle t_{2n},\sin k\tau\right\rangle=0,
  \end{array}\right.k=0,1,2\cdots,n-1,
\end{equation}
which in turn implies
\begin{equation}
Y_{1,1}(z)=t_{2n}(z)=\varpi_{2n}(z),
\end{equation}
which is the first kind monic OLPs defined by (1.5). Then the
solution of Riemann-Hilbert problem (2.5) must be rewritten as
\begin{equation}
Y_{1,2}(z)=\f{e^{-inz}}{4\pi
i}\int_0^{2\pi}\varpi_{2n}(\tau)\cot\f{\tau-z}2w(\tau)\mathrm{d}\tau,\,\,\,z\in
\mathbb{C}\setminus\mathbb{R},
\end{equation}
or say
\begin{equation}
Y_{1,2}(z)=\left\{
\begin{array}{ll}
 \displaystyle  \sum^\?_{k=n}\f1{2\pi}\!\!\int^{2\pi}_0t_{2n}(\tau)e^{-ik\tau}w(\tau)\mathrm{d}\tau e^{-i(n-k)z},\,\,\,\,&\mathrm{Im} z>0,
\\[3mm]
\displaystyle
-\sum^\?_{k=n}\f1{2\pi}\!\!\int^{2\pi}_0t_{2n}(\tau)e^{ik\tau}w(\tau)\mathrm{d}\tau
e^{-i(n+k)z},\,\,\,\,&\mathrm{Im} z<0.
\end{array}\right.
\end{equation}

On the contrary, reversing from (2.20) to (2.16), we easily know that (2.24) is exactly the unique solution of (2.5).

Secondly, we will solve Riemann-Hilbert problems (2.6) and (2.7). Similarly to the preceding discussion,
the solution of Riemann-Hilbert problem (2.6)
can be explicitly expressed by
\begin{equation}
  Y_{2,1}(z)=t_{2n-1}(z)\in T_{2n-1}(\mathbb{C}).
\end{equation}
If (2.7) is solvable, let $\widetilde{Y}_{2,2}(z)= e^{inz}Y_{2,2}(z)$. Then we have
\begin{equation}
\left\{
\begin{array}{l}
 \widetilde{Y}_{2,2}^+(x)= \widetilde{Y}_{2,2}^-(x)+t_{2n-1}(x)w(x),\,\,\,\,x\in [0,2\pi],
\\[3mm]
\displaystyle \widetilde{Y}_{2,2}(\pm \infty i)=0.
\end{array}\right.
\end{equation}
The unique solution of Riemann-Hilbert problem (2.27) can be written as
\begin{equation}
\widetilde{Y}_{2,2}(z)=\mathcal{C}_w[t_{2n-1}](z)=\f1{4\pi i}\int_0^{2\pi}t_{2n-1}(\tau)\cot\f{\tau-z}2w(\tau)\mathrm{d}\tau,\,\,\,z\in \mathbb{C}\setminus\mathbb{R},
\end{equation}
which in turn implies
\begin{equation}
Y_{2,2}(z)=\f{e^{-inz}}{4\pi i}\int_0^{2\pi}t_{2n-1}(\tau)\cot\f{\tau-z}2w(\tau)\mathrm{d}\tau,\,\,\,z\in \mathbb{C}\setminus\mathbb{R}.
\end{equation}
Further, putting (2.19) into (2.29), one has
\begin{equation}
Y_{2,2}(z)=\left\{
\begin{array}{l}
 \displaystyle  \f1{4\pi}\!\!\int^{2\pi}_0t_{2n-1}(\tau)w(\tau)\mathrm{d}\tau e^{-inz}+\sum^\?_{k=1}\f1{2\pi}\!\!\int^{2\pi}_0t_{2n-1}(\tau)e^{-ik\tau}w(\tau)\mathrm{d}\tau e^{-i(n-k)z},\,\,\,\,\mathrm{Im} z>0,
\\[3mm]
\displaystyle -\f1{4\pi}\!\!\int^{2\pi}_0t_{2n-1}(\tau)w(\tau)\mathrm{d}\tau e^{-inz}-\sum^\?_{k=1}\f1{2\pi}\!\!\int^{2\pi}_0t_{2n-1}(\tau)e^{ik\tau}w(\tau)\mathrm{d}\tau e^{-i(n+k)z},\,\,\,\,\mathrm{Im} z<0.
\end{array}\right.
\end{equation}
Thus, considering two growth conditions in (2.7), one easily knows, if and only if
\begin{equation}
\left\{
\begin{array}{l}
 \displaystyle \int^{2\pi}_0t_{2n-1}(\tau)e^{-ik\tau}w(\tau)\mathrm{d}\tau=0,\ k=0,1,2\cdots,n-2,
\\[3mm]
\displaystyle \f1{2\pi}\int^{2\pi}_0t_{2n-1}(\tau)e^{-i(n-1)\tau}w(\tau)\mathrm{d}\tau=1
\end{array}\right.
\end{equation}
and
\begin{equation}
\left\{
\begin{array}{l}
 \displaystyle \int^{2\pi}_0t_{2n-1}(\tau)e^{ik\tau}w(\tau)\mathrm{d}\tau=0,\ k=1,2\cdots,n-2,
\\[3mm]
\displaystyle -\f1{2\pi}\int^{2\pi}_0t_{2n-1}(\tau)e^{i(n-1)\tau}w(\tau)\mathrm{d}\tau=1,
\end{array}\right.
\end{equation}
the unique solution of (2.7) can be expressed by (2.29) or (2.30). Further, (2.31) and (2.32) are equivalent to
\begin{equation}
   \left\{\begin{array}{l}
   \displaystyle    \int^{2\pi}_0t_{2n-1}(\tau)\cos k\tau w(\tau)\mathrm{d}\tau=0,\ k=0,1,2\cdots,n-1,\\[2mm]
  \displaystyle    \int^{2\pi}_0t_{2n-1}(\tau)\sin k\tau w(\tau)\mathrm{d}\tau=0,\ k=1,2\cdots,n-2,\\[2mm]
  \displaystyle   \int^{2\pi}_0t_{2n-1}(\tau)\sin(n-1)\tau w(\tau)\mathrm{d}\tau=2\pi i,
    \end{array}
    \right.
 \end{equation}
which leads to
\begin{equation}
Y_{2,1}(z)=t_{2n-1}(z)=a_n\varpi_{2n-1}(z),
\end{equation}
 where $a_n$ is defined by (2.12).

 Finally, inserting (2.34) into (2.29), one gets
\begin{equation}
Y_{2,1}(z)=a_ne^{-inz}\mathcal{C}_w[\varpi_{2n-1}](z).
\end{equation}

Similarly, reversing step by step, we obtain that (2.35) is just the unique solution of (2.7). This completes the proof.\,\,\,\,\,\,\,$\blacksquare$

\section{the steepest descent analysis}

In this section, we will carry out an array of transforms
\begin{equation}
\mathbf{Y}\mapsto \mathbf{F}\mapsto \mathbf{S}\mapsto \mathbf{R},
\end{equation}
and the model Riemann-Hilbert problem is obtained. And the steepest descent analysis bases on those transforms.

\subsection{}the first transform $\mathbf{Y}\mapsto \mathbf{F}$
\vskip5mm

Let
\begin{equation}
  \Gamma(z)=\f1{4\pi i}\int^{2\pi}_0\ln w(\tau)\cot\f{\tau-z}2\mathrm{d}\tau,\ z\notin \mathbb{R}.
\end{equation}
Clearly,
\begin{equation}
  \Gamma(\pm \?i)=\pm\f1{4\pi }\int^{2\pi}_0\ln w(\tau)\mathrm{d}\tau.
\end{equation}

We define
\begin{equation}
   \left\{\begin{array}{ll}
    D^+(z)=e^{\Gamma(z)-C},\ &\mathrm{Im}z>0,\\[2mm]
    D^-(z)=e^{-\Gamma(z)-C},\ &\mathrm{Im}z<0
    \end{array}
    \right.
 \end{equation}
 with
 \begin{equation}\label{ccc}
  C=\f1{4\pi }\int^{2\pi}_0\ln w(\tau)\mathrm{d}\tau.
\end{equation}
Obviously,
\begin{equation}
   \left\{\begin{array}{l}
    D^+\in A(\square_{0,+\infty}),\,\,\,D^-\in A(\square_{-\infty,0}),\\[3mm]
    D^+(z+2\pi)= D^+(z),\,\,\,D^-(z+2\pi)= D^-(z),\\[3mm]
      D^+(x)D^-(x)=e^{-2C}w(x),\,\,\,x\in [0,2\pi],\\[3mm]
       D^\pm(\pm \infty i)=1,
    \end{array}
    \right.
 \end{equation}
 where
\begin{equation} \left\{\begin{array}{l}
  \square_{a,+\infty}=\big\{z,\mathrm{Re}z\in (0,2\pi),\,\mathrm{Im}z\in (a,+\infty)\big\}\\[3mm]
  \square_{-\infty,a}=\big\{z,\mathrm{Re}z\in (0,2\pi),\,\mathrm{Im}z\in (-\infty,a)\big\}
  \end{array}\right.\,\,\,\mathrm{with}\,\,\,a\in \mathbb{R}
\end{equation}
are infinite rectangle domains.

Further, we also define
\begin{equation}
    \mathfrak{D}^+(z)=\left\{\begin{array}{ll}
    D^+(z),\ &\mathrm{Im}z\geq 0,\\[2mm]
    e^{-2C}w(z)[D^-(z)]^{-1},\ &\mathrm{Im}z\in (-\rho,0],
    \end{array}
    \right.
 \end{equation}
 \begin{equation}
    \mathfrak{D}^-(z)=\left\{\begin{array}{ll}
    D^-(z),\ &\mathrm{Im}z\leq 0,\\[2mm]
    e^{-2C}w(z)[D^+(z)]^{-1},\ &\mathrm{Im}z\in [0,\rho),
    \end{array}
    \right.
 \end{equation}
 where $\rho$ is consistent with that in (1.12). And hence, by (1.12) and (3.6),
 \begin{equation}
   \left\{\begin{array}{l}
    \mathfrak{D}^+\in A(\square_{-\rho,+\infty}),\,\,\,\mathfrak{D}^-\in A(\square_{-\infty,\rho}),\\[3mm]
    \mathfrak{D}^+(z+2\pi)= \mathfrak{D}^+(z),\,\,\,\mathfrak{D}^-(z+2\pi)= \mathfrak{D}^-(z),\\[3mm]
      \mathfrak{D}^+(x)\mathfrak{D}^-(x)=e^{-2C}w(x),\,\,\,x\in [0,2\pi],\\[3mm]
       \mathfrak{D}^\pm(\pm \infty i)=1.
    \end{array}
    \right.
 \end{equation}

 Now, set
 \begin{equation}
  \mathbf{U}(z)=\left\{
  \begin{array}{ll}
\left(
           \begin{array}{cc}
             \f{e^{inz}\mathfrak{D}^+(z)}2 & 0 \\
             0 &  \f{e^{iz}}{\mathfrak{D}^+(z)} \\
           \end{array}
         \right),\ &z\in \square_{0,+\infty},\\[5mm]
    \left(
           \begin{array}{cc}
             \f{e^{-inz}\mathfrak{D}^-(z)}2 & 0 \\
             0 &  \f{e^{i(2n-1)z}}{\mathfrak{D}^-(z)} \\
           \end{array}
         \right)  ,\ &z\in \square_{-\infty,0},
         \end{array}\right.
\end{equation}
and we define the first transform
\begin{equation}
\mathbf{F}(z)=\mathbf{Y}(z)\mathbf{U}(z),\,\,\,z\in \mathbb{C}\setminus \mathbb{R}.
 \end{equation}

\noindent{\bf Lemma 3.1}\,\,\,{\it If $\mathbf{Y}$ is the solution of Riemann-Hilbert problem (2.1), then $\mathbf{F}(z)$ defined by (3.12) is the solution of the following Riemann-Hilbert problem
\begin{equation}
 \left\{\begin{array}{lll}
\mathbf{F}^+(x)=\mathbf{F}^-(x)\left(
                \begin{array}{cc}
                  e^{i2nx}\f{[\mathfrak{D}^+(x)]^2}{w(x)} & 2e^{ix}\\
                  0 & e^{-i2(n-1)x}\f{[\mathfrak{D}^-(x)]^{2}}{w(x)}\\
                \end{array}
              \right)e^{2C},\,\,\,&x\in [0,2\pi],\\[5mm]
  \mathbf{F}(z)=\mathbf{I}+\mathbf{o}(1),\ &z\rt+\?i,\\[5mm]
            \mathbf{F}(z)=\mathbf{I}+\mathbf{o}(1),\ &z\rt-\?i.
            \end{array}\right.
\end{equation}
Conversely, if $\mathbf{F}(z)$ is the solution of the Riemann-Hilbert problem (3.13), then
\begin{equation}
\mathbf{Y}(z)=\mathbf{F}(z)\mathbf{U}^{-1}(z),\,\,\,z\in \mathbb{C}\setminus \mathbb{R}
 \end{equation}
is the solution of Riemann-Hilbert problem (2.1).}\vskip5mm

\subsection{} The second transform $\mathbf{F}(z)\longrightarrow \mathbf{S}(z) $
\vskip5mm

The coefficient matrix in the boundary condition in (3.13) can be decomposed as follows
\begin{eqnarray}
  \nonumber &\ & \left(
                \begin{array}{cc}
                  e^{i2nx}\f{[\mathfrak{D}^+(x)]^2}{w(x)} & 2e^{ix} \\
                  0 & e^{-i2(n-1)x}\f{[\mathfrak{D}^-(x)]^{2}}{w(x)}\\
                \end{array}
              \right)\,\,\,\,\,(x\in [0,2\pi])
   \\
   &=& \left(
                 \begin{array}{cc}
                   1 & 0 \\
                   \f12e^{-i(2n-1)x}\f{[\mathfrak{D}^-(x)]^2}{w(x)} & 1 \\
                 \end{array}
               \right)\left(
                        \begin{array}{cc}
                          0 & 2e^{ix} \\
                          -\f12e^{ix} & 0 \\
                        \end{array}
                      \right)
               \left(
                               \begin{array}{cc}
                                 1 & 0 \\
                                 \f12e^{i(2n-1)x}\f{[\mathfrak{D}^+(x)]^2}{w(x)} & 1 \\
                               \end{array}
                             \right).
\end{eqnarray}
And hence the boundary condition in (3.13) is changed to
 \begin{equation}\label{128}
  \mathbf{F}^+(x)\left(
                               \begin{array}{cc}
                                 1 & 0 \\
                                 -\f12e^{i(2n-1)x}\f{[\mathfrak{D}^+(x)]^2}{w(x)} & 1 \\
                               \end{array}
                             \right)=\mathbf{F}^-(x)\left(
                 \begin{array}{cc}
                   1 & 0 \\
                   \f12e^{-i(2n-1)x}\f{[\mathfrak{D}^-(x)]^2}{w(x)} & 1 \\
                 \end{array}
               \right)\left(
                        \begin{array}{cc}
                          0 & 2e^{ix} \\
                          -\f12e^{ix} & 0 \\
                        \end{array}
                      \right)e^{2C}.
\end{equation}

In order to construct the second transform, we introduce some symbols. For $r\in (0,\rho)$, we define two oriented line segments and a contour
\begin{equation}
L_r=\left[2\pi+ir, ir\right],\,\,\,\,L_{-r}=\left[2\pi-ir, -ir\right],\,\,\,\Gamma=L_r+[0,2\pi]+L_{-r}.
\end{equation}
And the contour $\Gamma$ divides the basic strip
\begin{equation}
\square_{-\infty,+\infty}=\left\{z,\,\mathrm{Re}z\in (0,2\pi),\,\mathrm{Im}z\in (-\infty,+\infty)\right\}
\end{equation}
into two domains
\begin{equation}
\mathbb{A}^+=\mathbb{A}_1^++\mathbb{A}_2^+,\,\,\,\,\mathbb{A}^-=\mathbb{A}_1^-+\mathbb{A}_2^-
\end{equation}
with
\begin{equation}
\mathbb{A}_1^+=\square_{-\infty,-r},\,\,\mathbb{A}_2^+=\square_{0,r},\,\,\,\,\mathbb{A}_1^-=\square_{-r,0},\,\,\, \mathbb{A}_2^-=\square_{r,+\infty}.
\end{equation}

Let
\begin{equation}
  \mathbf{V}(z)=\left\{
  \begin{array}{ll}
  \mathbf{I},&z\in \mathbb{A}_1^+,\,\\[3mm]
  \left(
        \begin{array}{cc}
          1 & 0 \\
          \f12e^{-i(2n-1)z}\f{[\mathfrak{D}^-(z)]^2}{w(z)} & 1 \\
        \end{array}
      \right),\ &z\in \mathbb{A}_1^-,\\[5mm]
     e^{-iz}\left(
        \begin{array}{cc}
          1 & 0 \\
          -\f12e^{i(2n-1)z}\f{[\mathfrak{D}^+(z)]^2}{w(z)} & 1 \\
        \end{array}
      \right),\ &z\in \mathbb{A}_2^+,\\[3mm]
   e^{-iz}\mathbf{I}  &z\in \mathbb{A}_2^-,
      \end{array}\right.
\end{equation}
and we define the second transform
\begin{equation}
\mathbf{S}(z)=\mathbf{F}(z)\mathbf{V}(z),\,\,\,z\in \mathbb{C}\setminus \mathbb{R}.
 \end{equation}

Similarly to Lemma 3.1, one has the following. \vspace{3mm}

\noindent{\bf Lemma 3.2}\,\,\,{\it If $\mathbf{F}$ is the solution of Riemann-Hilbert problem (3.13), then $\mathbf{S}(z)$ defined by (3.22) is the following Riemann-Hilbert problem
\begin{equation}
  \left\{
  \begin{array}{ll}
  \mathbf{S}^+(t)=\mathbf{S}^-(t)\mathbf{\Upsilon}(t),\ &t\in \Gamma,\\[3mm]
  \mathbf{S}(z)=e^{-iz}\big[\mathbf{I}+\mathbf{o}(1)\big],\ &z\rt+\?i,\\[5mm]
           \mathbf{S}(z)=\mathbf{I}+\mathbf{o}(1),\ &z-\rt-\?i
                   \end{array}\right.
\end{equation}
with
\begin{equation}
\mathbf{\Upsilon}(t)=\left\{
  \begin{array}{ll}
  \left(
        \begin{array}{cc}
          1 & 0 \\
          \f12e^{i(2n-1)t}\f{[\mathfrak{D}^+(t)]^2}{w(t)} & 1 \\
        \end{array}
      \right),&t\in L_r,\\[5mm]
      \left(
        \begin{array}{cc}
          0 & 2 \\
          -\f12& 0 \\
        \end{array}
      \right)^{-1}e^{-2C},\ &t\in [0,2\pi],\\[5mm]
       \left(
        \begin{array}{cc}
          1 & 0 \\
         \f12 e^{-i(2n-1)t}\f{[\mathfrak{D}^-(t)]^2}{w(t)} & 1 \\
        \end{array}
      \right),&t\in L_{-r}.
      \end{array}\right.
\end{equation}
Conversely, if $\mathbf{S}(z)$ is the solution of the Riemann-Hilbert problem (3.23), then
\begin{equation}
\mathbf{F}(z)=\mathbf{S}(z)\mathbf{V}^{-1}(z),\,\,\,z\in \mathbb{C}\setminus \mathbb{R}
 \end{equation}
is the solution of Riemann-Hilbert problem (3.13).}\vskip5mm

\subsection{} The third transform $\mathbf{S}(z)\longrightarrow \mathbf{R}(z)$
\vskip5mm

To eliminate the jump of $\mathbf{S}(z)$ on the interval $[0,2\pi]$, we need to find the sectionally periodic analytic function $\mathbf{M}(z)$ satisfying the following conditions:
\begin{equation}\label{143}
  \left\{
  \begin{array}{ll}
  \mathbf{M}^+(x)=\mathbf{M}^-(x)\left(
                 \begin{array}{cc}
                   0 & 2 \\
                   -\f12 & 0 \\
                 \end{array}
               \right)e^{2C},\ & x\in [0,2\pi],\\[5mm]
              \mathbf{M}(z)=\left(
                 \begin{array}{cc}
                   0 & 2 \\
                   -\f12 & 0 \\
                 \end{array}
               \right)e^{2C}+\mathbf{o}(1),\ &z\rt +\?i,\\[4mm]
\mathbf{M}(z)=\mathbf{I}+\mathbf{o}(1),\ &z\rt-\?i,
                        \end{array}\right.
\end{equation}
where $C$ is defined by (3.5).  By Liouvile's Theorem,
 one easily knows that
\begin{equation}\label{1431}
 \mathbf{ M}(z)=\left\{
  \begin{array}{ll}

              \left(
                 \begin{array}{cc}
                   0 & 2 \\
                   -\f12 & 0 \\
                 \end{array}
               \right)e^{2C},&\ \mathrm{Im}z>0, \\[5mm]
             \left(
                     \begin{array}{cc}
                       1 & 0 \\
                       0 & 1 \\
                     \end{array}
                   \right),&\ \mathrm{Im}z<0
\end{array}
\right.
\end{equation}
is the unique solution of (3.26).

Now, the contour
\begin{equation}
\Gamma_\sharp=L_r+L_{-r}^-
\end{equation}
divides
the basic strip
$
\square_{-\infty,+\infty}
$
into two domains
\begin{equation}
\mathbb{B}^+=\square_{-r,r}\,\,\,\mathrm{and}\,\,\,\,\mathbb{B}^-=\mathbb{B}_1^-+\mathbb{B}_2^-
\end{equation}
with
\begin{equation}
\mathbb{B}_1^-=\square_{-\infty,-r},\,\,\,\, \mathbb{B}_2^-=\square_{r,+\infty}.
\end{equation}
We define the third transform
\begin{equation}
\mathbf{R}(z)=\mathbf{S}(z)\mathbf{M}^{-1}(z),\,\,\,z\in \mathbb{B}^+\cup \mathbb{B}^-,
\end{equation}
where $\mathbf{M}^{-1}(z)$ is the inverse of $\mathbf{M}(z)$ defined by (3.27).
By the boundary conditions in (3.23) and (3.26), one has
\begin{equation}
\mathbf{S}^+(x)\big[\mathbf{M}^+(x)\big]^{-1} =\mathbf{S}^-(x)\big[\mathbf{M}^-(x)\big]^{-1},\ x\in [0,2\pi],
\end{equation}
which implies that $\mathbf{R}(z)$ can be analytically extended across $\mathbb{R}$. And hence, we always assume that $\mathbf{R}(z)$ is analytic on $\mathbb{R}$ in what follows.

Further, by a simple calculation, one easily gets
\begin{eqnarray}
  \nonumber \mathbf{R}^-(t) &=& \mathbf{S}^-(t)\mathbf{M}^{-1}(t)=\mathbf{S}^+(t)\left(
                                      \begin{array}{cc}
                                        1 & 0 \\
                                        -\f12e^{i(2n-1)t}\f{[\mathfrak{D}^+(t)]^2}{w(t)} & 1\\
                                      \end{array}
                                    \right)\mathbf{M}^{-1}(t)
  \\
  \nonumber  &=& \mathbf{S}^+(t)\mathbf{M}^{-1}(t)\mathbf{M}(t)\left(
                                      \begin{array}{cc}
                                        1 & 0 \\
                                        -e^{i(2n-1)t}\f{[\mathfrak{D}^+(t)]^2}{w(t)} & 1\\
                                      \end{array}
                                    \right)\mathbf{M}^{-1}(t)
                                    \\
            \nonumber &=&    \mathbf{R}^+(t)                     \left(
                 \begin{array}{cc}
                   0 & 2 \\
                   -\f12 & 0 \\
                 \end{array}
               \right)e^{2C } \left(
                                      \begin{array}{cc}
                                        1 & 0 \\
                                        -\f12e^{i(2n-1)t}\f{[\mathfrak{D}^+(t)]^2}{w(t)} & 1\\
                                      \end{array}
                                    \right)\left(
                 \begin{array}{cc}
                   0 & -2 \\
                   \f12 & 0 \\
                 \end{array}
               \right)e^{-2C }\\
   &=& \mathbf{R}^+(t)\left(
               \begin{array}{cc}
                 1 & 2e^{i(2n-1)t}\f{[\mathfrak{D}^+(t)]^2}{w(t)}  \\
                 0 & 1 \\
               \end{array}
             \right),\,\,\,\,t\in L_r
\end{eqnarray}
and
\begin{eqnarray}\label{}
  \nonumber \mathbf{R}^+(t) &=&\mathbf{S}^+(t)\mathbf{M}^{-1}(t)  \\
  \nonumber  &=&\mathbf{S}^-(t)\left(
                     \begin{array}{cc}
                       1 & 0 \\
                      \f12 e^{-i(2n-1)t}\f{[\mathfrak{D}^-(t)]^2}{w(t)} & 1 \\
                     \end{array}
                   \right)\mathbf{M}^{-1}(t)
    \\
  \nonumber  &=&  \mathbf{S}^-(t)\mathbf{M}(t)\mathbf{M}^{-1}(t)\left(
                     \begin{array}{cc}
                       1 & 0 \\
                       \f12e^{-i(2n-1)t}\f{[\mathfrak{D}^-(t)]^2}{w(t)} & 1 \\
                     \end{array}
                   \right)\mathbf{M}^{-1}(t)\\
                   &=&\mathbf{R}^-(t)\left(
                     \begin{array}{cc}
                       1 & 0 \\
                      \f12 e^{-i(2n-1)t}\f{[\mathfrak{D}^-(t)]^2}{w(t)} & 1 \\
                     \end{array}
                   \right),\,\,\,t\in L_{-r}.
\end{eqnarray}
And then, the following lemma is obtained.
\vskip3mm

\noindent{\bf Lemma 3.3}\,\,\,{\it If $\mathbf{S}$ is the solution of Riemann-Hilbert problem (3.23), then $\mathbf{R}(z)$ defined by (3.31) is the solution of the model Riemann-Hilbert problem
\begin{equation}\label{171}
  \left\{
  \begin{array}{ll}
  \mathbf{R}^+(t)=\mathbf{R}^-(t)\mathbf{G}(t),& \ t\in \Gamma_\sharp=L_{r}+L_{-r}^-,\\[3mm]
  \mathbf{R}(z)=e^{-2C-iz }\left[\left(
                 \begin{array}{cc}
                   0 & -2 \\
                   \f12 & 0 \\
                 \end{array}
               \right)+\mathbf{o}(1)\right],\ &z\rt+\?i,\\[5mm]
               \mathbf{R}(z)=\mathbf{I}+\mathbf{o}(1),&z\rt-\?i.
                            \end{array}\right.
\end{equation}
with
\begin{equation}\label{162}
  \mathbf{G}(t)=\left\{
  \begin{array}{ll}
  \left(
    \begin{array}{cc}
      1 & -2e^{i(2n-1)t}\f{[\mathfrak{D}^+(t)]^2}{w(t)} \\
      0 & 1 \\
    \end{array}
  \right),&\ t\in L_{r},\\[5mm]
  \left(
    \begin{array}{cc}
      1 & 0 \\
      -\f12e^{-i(2n-1)t}\f{[\mathfrak{D}^-(t)]^2}{w(t)}  & 1 \\
    \end{array}
  \right),&t\in L_{-r}^{-}.
  \end{array}\right.
\end{equation}
Conversely, if $\mathbf{R}(z)$ is the solution of the model Riemann-Hilbert problem (3.35), then
\begin{equation}
\mathbf{S}(z)=\mathbf{R}(z)\mathbf{M}(z),\,\,\,z\in \mathbb{C}\setminus \Gamma_\sharp
 \end{equation}
is the solution of Riemann-Hilbert problem (3.23).}\vskip5mm

\section{Strong Asymptotic analysis of OTP}

First, we set up a lemma needed in the sequel. \vskip3mm

\bL Let $\mathbf{R}^-$ is the negative boundary value of
$\mathbf{R}$ which is the solution for the model Riemann-Hilbert
problem (3.35). Then
\begin{equation}\label{}
  \mathbf{k}=\f1{4\pi }\int_{\Gamma_\sharp} \mathbf{R}^-(\tau)(\mathbf{G}(\tau)-\mathbf{I})\mathrm{d}\tau=-\frac 12 \mathbf{I},
          \end{equation}
where $\mathbf{G}$ is given by (3.36) and $\mathbf{I}$ is the identity matrix.
\eL

\bpp Recall that $\mathbf{R}^+$ is analytic on $ \square_{-r,r}$. We have
\begin{eqnarray}\label{lambda}
\nonumber  \f1{4\pi }\int_{\Gamma_\sharp}\mathbf{R}^-(\tau)(\mathbf{G}(\tau)-\mathbf{I})\mathrm{d}\tau
&=& \f1{4\pi }\int_{\Gamma_\sharp} \mathbf{R}^-(\tau)\mathbf{G}(\tau)\mathrm{d}\tau
-\f1{4\pi }\int_{\Gamma_\sharp} \mathbf{R}^-(\tau)\mathrm{d}\tau \\
 \nonumber  &=&  \f1{4\pi }\int_{\partial \square_{-r,r}} \mathbf{R}^+(\tau)\mathrm{d}\tau
 - \f1{4\pi }\int_{\Gamma_\sharp} \mathbf{\mathbf{R}}^-(\tau)\mathrm{d}\tau\\
  &=& -\f1{4\pi }\int_{L_{-r}} \mathbf{R}^-(\tau)\mathrm{d}\tau+\f1{4\pi }\int_{L_{r}}
  \mathbf{R}^-(\tau)\mathrm{d}\tau,
\end{eqnarray}
where $L_r,\,L_{-r}$ are defined in (3.17). But by (\ref{171}) we
have
\begin{eqnarray}\label{lambda2}
  \nonumber\f1{4\pi }\int_{L_{-r}}
  \mathbf{R}^-(\tau)\mathrm{d}\tau &=& \lim_{R\rt+\?} \f1{4\pi }\int_{L_{-R}} \mathbf{R}(z)\mathrm{d}z\\
\nonumber   &=& \lim_{R\rt+\?} \f1{4\pi }\int_{L_{- R}}
\left(\mathbf{I}
+\mathbf{o}(1)\right)\mathrm{d}z \\
   &=& \f12\, \mathbf{I}.
\end{eqnarray}

Denote by $A(z)=\mathbf{R}(z)e^{2C+iz }$ and $B(w)=A(-i\ln w)$ for $w\in \CC$. Then $B(w)$ is analytic near $\infty$. And  by (\ref{171}) its Laurant series is
 $$
 B(w)=\left(
                 \begin{array}{cc}
                   0 & -2 \\
                   \f12 & 0 \\
                 \end{array}
               \right)+\f1w c_1+\f1{w^2}c_2+\ldots,
               $$
 where $c_j,j=1,2,\cdots$ are constant matrices. So we have
 \begin{equation}\label{matrix}
   \mathbf{R}(z)=e^{-2C-iz }\left[\left(
                 \begin{array}{cc}
                   0 & -2 \\
                   \f12 & 0 \\
                 \end{array}
               \right)+c_1e^{iz}+\mathbf{o}(|e^{iz}|)\right].
 \end{equation}
 Then we obtain
 \begin{eqnarray}\label{lambda1}
 \nonumber \f1{4\pi }\int_{L_{r}}
  \mathbf{R}^-(\tau)\mathrm{d}\tau &=& \lim_{R\rt+\?} \f1{4\pi }\int_{L_{R}} \mathbf{R}(z)\mathrm{d}z\\
 \nonumber  &=& \lim_{R\rt+\?} \f1{4\pi }\int_{L_{R}} e^{-2C-iz }\left[\left(
                 \begin{array}{cc}
                   0 & -2 \\
                   \f12 & 0 \\
                 \end{array}
               \right)+c_1e^{iz}+\mathbf{o}(|e^{iz}|\right]\mathrm{d}z.
\end{eqnarray}\
Finally, combining (\ref{lambda}), (\ref{lambda1}) with
(\ref{lambda2}), we get the conclusion of the lemma. \epp

Now, one comes to verify one of the main results, usually called Aptekarev type theorem \cite{Z-J1,Z-J2}. It must be pointed that all the norms in the following have the same definition with those
in \cite{Z-J1,Z-J2,W-L-D}.

\bT \label{main} There exist constants $\eta>0$ and $\delta>0$ such
that, for $\|\mathbf{G}-\mathbf{I}\|_{\Omega_\epsilon}<\delta$, we
have
\begin{equation}\label{theorem}
 \left\|\mathbf{R}(z)-e^{-2C-iz}\left(
                 \begin{array}{cc}
                   0 & -2 \\
                   \f12 & 0 \\
                 \end{array}
               \right)-\frac 12 \mathbf{I}\right\|_{\CC\setminus \Gamma_\sharp}
               < \eta\|\mathbf{G}-\mathbf{I}\|_{\Omega_\epsilon}\,\,\,\mathrm{with}\,\,\,\Omega_\epsilon=\square_{r-\epsilon,r+\epsilon}\cup \square_{-r-\epsilon,-r+\epsilon}
\end{equation}
where $\epsilon>0$ is sufficiently small, $\mathbf{G}$ is given by (3.36) and $\mathbf{R}$ is the solution for the model Riemann-Hilbert
problem (3.35).\eT

\bpp First, let
$$
\mathbf{W}(z)=\mathbf{R}(z)-e^{-2C-iz}\left(
                 \begin{array}{cc}
                   0 & -2 \\
                   \f12 & 0 \\
                 \end{array}
               \right),\ z\in \CC\setminus \Gamma_\sharp.
               $$
By the asymptotic conditions in (2.35), one has
\begin{equation}\label{1711}
  \mathbf{W}(z)=\left\{\begin{array}{ll}
  \mathbf{A}+\mathbf{o}(1),\ &z\rt+\?i,\\[5mm]
  \mathbf{I}+\mathbf{o}(1),&z\rt-\?i.
  \end{array}\right.
\end{equation}
Denote by $\mathbf{\Delta}=\mathbf{G}-\mathbf{I}$. Again by the
boundary condition in (\ref{171}), we have
\begin{equation}\label{delta}
  \mathbf{R}^+(t)=\mathbf{R}^-(t)+\mathbf{R}^-(t)\mathbf{\Delta}
  =\mathbf{R}^-(t)(\mathbf{I}+\mathbf{\Delta}) ,\ t\in \Gamma_\sharp,
\end{equation}
which leads to
\begin{equation}\label{01}
  \mathbf{W}^+(t)=\mathbf{W}^-(t)+\mathbf{R}^-(t)\mathbf{\Delta},\ t\in \Gamma_\sharp.
\end{equation}

Secondly, we set
\begin{equation}\label{04}
\mathbf{H}(z)=\f1{4\pi
i}\int_{\Gamma_\sharp}\mathbf{R}^-(\tau)\mathbf{\Delta}(\tau)\cot\f{\tau-z}2\mathrm{d}\tau,\
z\notin \Gamma_\sharp.
\end{equation}
By Lemma 4.1, one has
\begin{equation}\label{1713}
  \left\{
  \begin{array}{ll}
  \mathbf{H}^+(t)=\mathbf{H}^-(t)+\mathbf{R}^-(t)\mathbf{\Delta},& \ t\in \Gamma_\sharp\\[3mm]
  \mathbf{H}(z)=\mathbf{k}+\mathbf{o}(1),\ &z\rt+\?i,\\[3mm]
               \mathbf{H}(z)=-\mathbf{k}+\mathbf{o}(1),&z\rt-\?i,
                  \end{array}\right.
\end{equation}
where $\mathbf{k}$ is given by (4.1). By (\ref{01}) and (\ref{1713}),
we have $\mathbf{W}(z)-\mathbf{H}(z)$ is an entire function. And by
(\ref{1711}) and (\ref{1713}), $\mathbf{W}(z)-\mathbf{H}(z)$ is
bounded on the whole complex plane. Then, by Liouville Theorem, we
have
\begin{equation}\label{03}
  \mathbf{W}(z)-\mathbf{H}(z)\equiv \mathbf{I}+\mathbf{k}=\frac 12 \mathbf{I},\ z\in\CC,
\end{equation}
which leads to
\begin{equation}\label{}
\mathbf{W}(z)=\frac 12\mathbf{I}+\f1{4\pi
i}\int_{\Gamma_\sharp}\mathbf{R}^-(\tau)\mathbf{\Delta}(\tau)\cot\f{\tau-z}2\mathrm{d}\tau,\
z\notin \Gamma_\sharp.
\end{equation}
Therefore, one has
\begin{equation}\label{}
\mathbf{R}(z)=e^{-2C-iz}\left(
                 \begin{array}{cc}
                   0 & -2 \\
                   \f12 & 0 \\
                 \end{array}
               \right)+\frac 12\mathbf{I}+\f1{4\pi
i}\int_{\Gamma_\sharp}\mathbf{R}^-(\tau)\mathbf{\Delta}(\tau)\cot\f{\tau-z}2\mathrm{d}\tau,\
z\notin \Gamma_\sharp.
\end{equation}

Thirdly, for $\epsilon>0$, let
$$\Gamma_\sharp^\epsilon=L_{r+\frac \epsilon 2}+L_{-r-\frac \epsilon 2}.$$
By (4.13), one has
\begin{equation}\label{}
\mathbf{R}(z)=e^{-2C-iz}\left(
                 \begin{array}{cc}
                   0 & -2 \\
                   \f12 & 0 \\
                 \end{array}
               \right)+\frac 12\mathbf{I}+\f1{4\pi
i}\int_{\Gamma_\sharp^\epsilon}\mathbf{R}^-(\tau)\mathbf{\Delta}(\tau)\cot\f{\tau-z}2\mathrm{d}\tau,\
z\in \mathbb{B}^+=\square_{-r,r},
\end{equation}
which in particular implies
\begin{equation}\label{}
\mathbf{R}^+(t)=e^{-2C-it}\left(
                 \begin{array}{cc}
                   0 & -2 \\
                   \f12 & 0 \\
                 \end{array}
               \right)+\frac 12\mathbf{I}+\f1{4\pi
i}\int_{\Gamma_\sharp^\epsilon}\mathbf{R}^-(\tau)\mathbf{\Delta}(\tau)\cot\f{\tau-t}2\mathrm{d}\tau,\
t\in \Gamma_\sharp.
\end{equation}
Further, by (4.7) and (4.15), we get
\begin{equation}\label{}
\mathbf{R}^-(t)=\left[e^{-2C-it}\left(
                 \begin{array}{cc}
                   0 & -2 \\
                   \f12 & 0 \\
                 \end{array}
               \right)+\frac 12\mathbf{I}+\f1{4\pi
i}\int_{\Gamma_\sharp^\epsilon}\mathbf{R}^-(\tau)\mathbf{\Delta}(\tau)\cot\f{\tau-t}2\mathrm{d}\tau\right]\left[\mathbf{I}+\mathbf{\Delta}(t)\right]^{-1},\
t\in \Gamma_\sharp.
\end{equation}
This leads to the estimate
\begin{equation}\label{07}
  \|\mathbf{R}^-\|_{\Gamma_\sharp}\<\left(d_\epsilon+K_\epsilon\|\mathbf{\Delta}\|_{\Omega_\epsilon} \|\mathbf{R}^-\|_{\Gamma_\sharp^\epsilon}\right)\left\|(\mathbf{I}+\mathbf{\Delta})^{-1}\right\|_{\Omega_\epsilon}
\end{equation}
with
\begin{equation}d_\epsilon=\left\|e^{-2C-it}\left(
                 \begin{array}{cc}
                   0 & -2 \\
                   \f12 & 0 \\
                 \end{array}
               \right)+\frac12\mathbf{I}\right\|_{\Omega_\epsilon},\,\,\,K_\epsilon=\max_{t\in \Gamma_\sharp,\t\in \Gamma_\sharp^\epsilon}\left|\cot\f{\t-t}2\right|.
\end{equation}
By the Maximal Principle Theorem, we have
\begin{equation}
\|\mathbf{R}^-\|_{\Gamma_\sharp}\geq\|\mathbf{R}^-\|_{\Gamma_\sharp^\epsilon}.
\end{equation}
Notice that
\begin{equation}
\left\|(\mathbf{I}+\mathbf{\Delta})^{-1}\right\|_{\Omega_\epsilon}\<\sum_{n=0}^\?\|\mathbf{\Delta}\|^n_{\Omega_\epsilon}=\f1{1-\|\mathbf{\Delta}\|_{\Omega_\epsilon}}
\end{equation}
for $\|\mathbf{\Delta}\|_{\Omega_\epsilon}<1$. Combining (4.17), (4.19) with (4.20), one has
  \begin{equation}\label{07}
  \|\mathbf{R}^-\|_{\Gamma_\sharp}\<\frac{d_\epsilon+K_\epsilon\|\mathbf{\Delta}\|_{\Omega_\epsilon} \|\mathbf{R}^-\|_{\Gamma_\sharp}}{1-\|\mathbf{\Delta}\|_{\Omega_\epsilon}}
\end{equation}
for $\|\mathbf{\Delta}\|_{\Omega_\epsilon}<1$.
 Now, we pick $\delta\in (0,\frac 12)$ sufficiently small such that
\begin{equation}
\|\mathbf{\Delta}\|_{\Omega_\epsilon}<\delta\Longrightarrow\frac{K_\epsilon\|\mathbf{\Delta}\|_{\Omega_\epsilon} }{1-\|\mathbf{\Delta}\|_{\Omega_\epsilon}}\leq\frac 12.
\end{equation}
Therefore, combing (4.21) with (4.22), one has
\begin{equation}
\|\mathbf{\Delta}\|_{\Omega_\epsilon}<\delta\Longrightarrow\|\mathbf{R}^-\|_{\Gamma_\sharp}\<\frac{\frac{d_\epsilon }{1-\|\mathbf{\Delta}\|_{\Omega_\epsilon}}}{1-\frac{K_\epsilon\|\mathbf{\Delta}\|_{\Omega_\epsilon} }{1-\|\mathbf{\Delta}\|_{\Omega_\epsilon}}}\leq \frac{2d_\epsilon}{1-\delta}<4d_\epsilon.
\end{equation}

Finally, according to \cite{mp}, the singular integral operator is $H^\mu-$bounded. And hence there exists $M>0$ such that
\begin{equation}\label{mp}
 \left\| \int_{\Gamma_\sharp}  \mathbf{R}^-(\t)\mathbf{\Delta}(\tau)\cot\f{\t-t}2\mathrm{d}\t\right\|_{\Gamma_\sharp}\leq M K_\epsilon\|\mathbf{\Delta}\|_{\Omega_\epsilon}\|\mathbf{R}^-\|_{\Gamma_\sharp}
\end{equation}
since the integrand is differentiable. Combining (4.23) with (4.24), one easily gets the desired estimate (4.5).

 \epp

 \vskip10mm

\noindent\bf Acknowledgements \rm This paper is supported
by  NSFC11701597; NSFC11471250.

\end{document}